\documentclass[11pt,reqno]{amsart}
\usepackage{latexsym}
\usepackage{amsfonts}
\usepackage{amssymb}
\usepackage{amsmath}
\usepackage{amsthm}
\usepackage{eucal}
\usepackage{amscd}
\usepackage{mathrsfs}
\usepackage{graphics}
\usepackage{epsfig}
\usepackage{fancyhdr}
\usepackage{verbatim}
\usepackage{setspace}
\usepackage{color}
\usepackage{enumerate}
\usepackage{hyperref}
\usepackage{enumitem,mathrsfs}
\usepackage{soul}
\usepackage[normalem]{ulem}
\usepackage{multirow}
\usepackage{cancel}

\usepackage{tikz}
\usepackage{tikz-cd}

\usepackage{mathtools}

\usepackage{enumitem}
\setlist[enumerate]{
wide, 
nosep, 
labelwidth=*,
labelindent=0cm,
leftmargin=0cm}

\newtheorem{theorem}{Theorem}
\newtheorem{proposition}[theorem]{Proposition}

\newtheorem{lemma}[theorem]{Lemma}

\newcommand{\R}{\mathbb{R}}
\newcommand{\Q}{\mathbb{Q}}
\newcommand{\K}{\mathbb{K}}
\newcommand{\N}{\mathbb{N}}
\newcommand{\C}{\mathbb{C}}

\usepackage{graphicx}

\makeatletter
\newcommand*\bigcdot{\mathpalette\bigcdot@{.8}}
\newcommand*\bigcdot@[2]{\mathbin{\vcenter{\hbox{\scalebox{#2}{$\m@th#1\bullet$}}}}}
\makeatother

\title[Lotka--Volterra systems with high degree invariant algebraic curves]
{A family of planar Lotka--Volterra systems with invariant algebraic curves of arbitrary degree}

\author[J. Coyo-Guarachi]{}
\author[S. Rebollo-Perdomo]{}


\thanks{The authors were supported by Universidad del B\'{\i}o-B\'{\i}o Grant RE2320122}

\begin{document}

\setlength{\abovedisplayskip}{3pt}
\setlength{\belowdisplayskip}{3pt}

\maketitle

\medskip

\centerline{\scshape Javier Coyo-Guarachi$^{a}$,
\ \ 
Salom\'on Rebollo-Perdomo$^{b}$ }
\medskip
{\footnotesize
\centerline{$^{a}$Departamento de Matemática, Universidad de Tarapacá, Arica, Chile.}
\centerline{javier.coyo.guarachi@alumnos.uta.cl}
\centerline{$^{b}$Departamento de Matemática, Universidad del B\'io-B\'io, Concepción, Chile.}
\centerline{srebollo@ubiobio.cl}
} 

\medskip

\begin{abstract}
We introduce a new family of planar  Lotka--Volterra systems admitting explicit invariant algebraic curves of arbitrarily high degree. 
\end{abstract}

\section{Introduction and statement of results}
\label{Seccion-introduccion}
Darboux showed in 1878 that a planar polynomial differential system of degree $d$ with at least $d(d+1)/2+1$
invariant algebraic curves has a Darboux first integral~\cite{Da1,Da2}. This result was improved by Jouanolou in 1979, 
who proved that a planar polynomial differential system of degree $d$ with at least $d(d+1)/2+2$
invariant algebraic curves has a rational first integral \cite{J}. 
Since then, the detection and computation of invariant algebraic curves for planar polynomial differential systems 
has been a subject of intensive research; see, for instance, \cite{CLPZ,LZa,LZb} and the references therein. 
However, determining whether a concrete planar polynomial differential system has invariant algebraic curves 
and the possible degrees of such curves, can be difficult problems.

\smallskip
Motivated by these issues, Brunella and Mendes \cite{BM} and Lins-Neto \cite{LN} recalled, 
at the beginning of this century, a classical problem that goes back to Poincar\'e \cite{Po}: 
given $d \geq 2$, does there exist a positive integer $M(d)$ such that
if a polynomial differential system of degree $d$ has an invariant algebraic curve of
degree $\geq M(d)$, then it has a rational first integral? 

\smallskip
Several examples provide a negative answer to this question in the case $d=2$; see, for example, \cite{CrLl, LlVa} and 
references therein. In particular, in 2001, Moulin-Ollagnier constructed a countable family of three-dimensional Lotka--Volterra 
systems, each of which has an associated planar Lotka--Volterra system, such that each system possesses an irreducible algebraic 
solution whose degree cannot be uniformly bounded \cite{MO}. However, the explicit expressions for these algebraic solutions were not provided.

\smallskip
In this work, we consider the non-countable family of planar Lotka--Volterra differential systems: 
\begin{equation}
\label{eq:sist-LV}
		\dot{x}=x(1-x),\quad
		\dot{y}=y\left(n +bx-y\right),
\end{equation}
where $(x,y)\in \K^2$, with $\K=\C$ or $\R$, $n\in \N$, and $b\in \K$. 
We will show that this family of Lotka--Volterra systems, distinct from the 
one studied in \cite{MO}, admits explicit invariant algebraic curves of arbitrarily high degree. 
More precisely, we establish the following results.

\begin{proposition}
\label{prop-curva de grado n}
For $n \in \N$ and $b\in \K$, consider the polynomial of degree $n$:
\begin{equation}
\label{equ: curv_inv_fn}
F(x,y)=y(n-1)!
\sum_{\nu=0}^{n-1}(-1)^{n+\nu}
(n+b-\nu+1)_{\nu}\frac{x^{\nu}}{\nu!}+(b+1)_nx^n \in \K[x,y],
\end{equation}
where
$(n+b-\nu+1)_{\nu}$ is the Pochhammer symbol of $n+b-\nu+1$ and $\nu.$  
The algebraic curve $\{F=0\}$ is an invariant algebraic curve of degree $n$ of the Lotka--Volterra 
system~\eqref{eq:sist-LV}, with cofactor $K=n-nx-y$. 
\end{proposition}

\begin{theorem}
\label{teo-integrable}
	Each Lotka--Volterra system \eqref{eq:sist-LV}
	is Darboux integrable. Moreover, if $b\in \Q$, then system \eqref{eq:sist-LV}
	has a rational first integral and if $b\not \in \Q$, then system \eqref{eq:sist-LV}
	has not a rational first integral.
\end{theorem}

\section{Basic definitions and proofs}
A {\it planar polynomial differential system of degree $d$} on $\K^2$, with $\K=\C$ or $\R$, 
is a differential system of the form
\begin{equation*}
\dot{x}:=\frac{dx}{dt} = P(x, y), \quad  
\dot{y}:=\frac{dy}{dt} = Q(x, y),
\end{equation*}
where $(x,y)\in \K^2$ are the dependent variables, $t\in \K$ is the independent 
variable (the time), $P$ and $Q$ are polynomials in the variables $x$ and $y$, {\it i.e.}, 
$P,Q\in \K[x,y]$, and $d = \max \{\deg P, \deg Q\}$.
An {\it invariant algebraic curve of degree $n$} of this planar polynomial differential system is the zero-locus
$$
\{F=0\}:=\{(x,y)\in \K^2 \, | \, F(x,y)=0\}
$$
of a polynomial $F\in \K[x,y]$ of degree $n$, such that the polynomial $F$ satisfies the linear
partial differential equation
\begin{equation}
\label{invcur}
PF_x+QF_y=KF,
\end{equation}
for some $K\in \K[x,y]$. Here $F_x$ and $F_y$ denote the partial derivatives of
$F$ respect to $x$ and $y$, respectively. The polynomial $K$ is called the {\it cofactor} of $\{F=0\}$.

\smallskip
Recall  that \emph{the Pochhammer symbol, $(c)_n$, of $c \in \K$ and $n\in\N$} is defined as
$$ 
(c)_n:=c(c+1)(c+2)\cdots(c+n-1) ;\hspace{0.5cm}(c)_0:=1.
$$ 
In order to prove Proposition~\ref{prop-curva de grado n}, 
the following results will be useful.

\begin{lemma}
\label{lem-Pochhammer-ident}
The Pochhammer symbol satisfies the following property
$$
\frac{(c+1)_{\nu}}{(\nu-1)!}
+\frac{(c)_{\nu+1}}{\nu !}
=\frac{(c+\nu)\,(c+1)_{\nu}}{\nu!}.
$$
\end{lemma}
\begin{proof}
It follows from the definition of $(c)_n$ and straightforward computations.
\end{proof}

\begin{lemma}
\label{lem-first-identity}
Consider $F(x,y)$ as in equation \eqref{equ: curv_inv_fn}. Then, the following identity holds.
$$
(1-x)F_x-(n-nx-y)(b+1)_nx^{n-1}=-(n+b)\big(F-(b+1)_nx^n\big).
$$
\end{lemma}
\begin{proof}
From the definition of $F(x,y)$, a direct computation yields
\begin{multline*}
(1-x)F_x=-y(n-1)!(-1)^{n}(n+b)\\
\qquad \qquad -y(n-1)!\sum_{\nu=1}^{n-2}(-1)^{n+\nu}\bigg( 
		\frac{(n+b-\nu)_{\nu+1}}{\nu !}
		+\frac{(n+b-\nu+1)_{\nu}}{(\nu-1)!}
		\bigg)x^{\nu}\\
		\!\!+y(n-1)!(b+2)_{n-1}\frac{x^{n-1}}{(n-2)!}+(b+1)_n(n-nx)x^{n-1}.
\end{multline*}
By using Lemma~\ref{lem-Pochhammer-ident}, with $c=n+b-\nu$, previous equation becomes
\begin{multline*}
(1-x)F_x=-y(n-1)!(-1)^{n}(n+b)\\-
y(n-1)!\sum_{\nu=1}^{n-2}(-1)^{n+\nu}\bigg( 
		\frac{(n+b)(n+b-\nu+1)_{\nu}}{\nu !}
		\bigg)x^{\nu}\\
		+y(n-1)!(b+2)_{n-1}\frac{x^{n-1}}{(n-2)!}+(b+1)_n(n-nx)x^{n-1},
\end{multline*}
which can be written as
\begin{multline*}
(1-x)F_x=-(n+b)y(n-1)!\sum_{\nu=0}^{n-2}(-1)^{n+\nu}
		(n+b-\nu+1)_{\nu}
		\frac{x^{\nu}}{\nu !}\\
		+y(n-1)!(b+2)_{n-1}\frac{x^{n-1}}{(n-2)!}+(b+1)_n(n-nx)x^{n-1}.
\end{multline*}
We now focus on the last two terms. The combination
$$
y(n-1)!(b+2)_{n-1}\frac{x^{n-1}}{(n-2)!}+(b+1)_n(n-nx)x^{n-1}-(n-nx-y)(b+1)_nx^{n-1}
$$
simplifies to
$$
y(n-1)!\left(\frac{(b+2)_{n-1}}{(n-2)!}+\frac{(b+1)_{n}}{(n-1)!} \right)x^{n-1},
$$
which, by using Lemma~\ref{lem-Pochhammer-ident}, with $c=b+1$ and $\nu=n-1$, becomes
$$
y(n-1)!\left(\frac{(n+b)(b+2)_{n-1}}{(n-1)!}\right)x^{n-1}.
$$
Consequently, the expression $(1-x)F_x-(n-nx-y)(b+1)_nx^{n-1}$ reduces to
$$
-(n+b)y(n-1)!\sum_{\nu=0}^{n-2}(-1)^{n+\nu}
		(n+b-\nu+1)_{\nu}
		\frac{x^{\nu}}{\nu !}+
		y(n-1)!\left(\frac{(n+b)(b+2)_{n-1}}{(n-1)!}\right)x^{n-1},
$$
which can be written as
$$
-(n+b)y(n-1)!\sum_{\nu=0}^{n-1}(-1)^{n+\nu}
		(n+b-\nu+1)_{\nu}
		\frac{x^{\nu}}{\nu !}.
$$
By recalling  equation \eqref{equ: curv_inv_fn}, this last equation is
\begin{equation*}
-(n+b)\Big(F-(b+1)_nx^n\Big).
\end{equation*}
This concludes the proof.
\end{proof}

\begin{proof}[Proof of Proposition~\ref{prop-curva de grado n}]
From the definition of $F(x,y)$ in equation \eqref{equ: curv_inv_fn}, we have
$$ 
yF_y=F-(1+b)_nx^n. 
$$
Thus, by using the expression of system \eqref{eq:sist-LV}, the left-hand side of \eqref{invcur} becomes
\begin{equation}
\label{eq-fund}
x(1-x)F_x+y(n+bx-y)F_y=x(1-x)F_x+(n+bx-y)\big(F-(b+1)_nx^n\big)
\end{equation}
We can decompose the term $n+bx-y$ as
$$
n+bx-y=n-nx-y+nx+bx=(n-nx-y)+(n+b)x.
$$
Substituting this identity into the right-hand side of \eqref{eq-fund} gives
$$
x(1-x)F_x+\big((n-nx-y)+(n+b)x\big) \big(F-(b+1)_nx^n\big),
$$
or equivalently,
\begin{multline*}
x\Big((1-x)F_x-(n-nx-y)(b+1)_nx^{n-1}+(n+b)\big(F-(b+1)_nx^{n}\big)\Big)+(n-nx-y)F.
\end{multline*}
Therefore, to conclude the proof, it suffices to prove that
$$
(1-x)F_x-(n-nx-y)(b+1)_nx^{n-1}=-(n+b)\big(F-(b+1)_nx^{n}\big).
$$
By Lemma~\ref{lem-first-identity}, this identity holds, which completes the proof.
\end{proof}

\begin{proof}[Proof of Theorem~\ref{teo-integrable}]
By construction, $\{f_1=0\}$ and $\{f_2=0\}$,  with $f_1=x$ and $f_2=y$, are invariant algebraic curves, 
whose cofactors are $K_1=1-x$ and $K_2=n+Bx-y$, respectively. It is clear that $\{f_3=0\}$, with $f_3=1-x$, 
is also an invariant algebraic curve of the system, and a direct computation shows that its cofactor is $K_3=-x$. 
In addition, by Proposition~\ref{prop-curva de grado n}, we know that  $\{f_4=0\}$, with
$$
f_4=\{y(n-1)!\sum_{\nu=0}^{n-1}(-1)^{n+\nu}(n+b-\nu+1)_{\nu}\frac{x^{\nu}}{\nu!}+(b+1)_nx^n,
$$
is an invariant algebraic curve of the system, whose cofactor is $K_4=n-nx-y$. 
By Darboux theorem \cite[Theorem 2]{LZa}, the function 
$H=f_1^{\lambda_1}f_2^{\lambda_2}f_3^{\lambda_3}f_4^{\lambda_4}$ 
is a Darboux first integral of the system if there are $\lambda_i \in \K$, 
with $i=1,2,3,4$ no all of them zero, such that  
\begin{equation*}
		\lambda_1K_1+\lambda_2K_2+\lambda_3K_3+\lambda_4K_4=0,
\end{equation*}
which is equivalent to
\begin{equation*}
	(\lambda_1+n\lambda_2+n\lambda_4)+(-\lambda_1+b\lambda_2-\lambda_3-n\lambda_4)x+(-\lambda_2-\lambda_4)y=0.
\end{equation*}
The solution to this equation is
$$
	\lambda_1=0,\quad
	\lambda_3=(n+b)\lambda_2,\quad
	\lambda_4=-\lambda_2.
$$
Hence, a Darboux first integral of the Lotka--Volterra system is
$$ 
H(x,y)=\frac{y^{\lambda_2}(1-x)^{(n+b)\lambda_2}}{f_4^{\lambda_2}}. 
$$
On the one hand, if $b=p/q\in \Q$, then for $\lambda_2=q$, 
$H(x,y)$ becomes a rational function. On the other hand, if $b \not \in \Q$, then $H(x,y)$ 
is not a rational function for any $\lambda_2\in \K$.
\end{proof}	

An analogous analysis can be used to prove the following two results.

\begin{proposition}
\label{prop-curva de grado n-2}
For $n \in \N$ and $b\in \K$, consider the polynomial of degree $n$:
\begin{equation*}
F(x,y)=-y(n-1)!
\sum_{\nu=0}^{n-1}(-1)^{n+\nu}
(n+b-\nu+1)_{\nu}\frac{x^{\nu}}{\nu!}+(b+1)_nx^n \in \K[x,y],
\end{equation*}
where
$(n+b-\nu+1)_{\nu}$ is the Pochhammer symbol of $n+b-\nu+1$ and $\nu.$  The  algebraic curve $\{F=0\}$ is an invariant algebraic curve of degree $n$ of the Lotka--Volterra system 
\begin{equation}
\label{eq:sist-LV-2}
		\dot{x}=x(1-x),\quad
		\dot{y}=y\left(n+bx+y\right),
\end{equation}
 with cofactor $K=n-nx+y$. 
\end{proposition}

\begin{theorem}
\label{teo-integrable-2}
	The Lotka--Volterra system \eqref{eq:sist-LV-2}
	is Darboux integrable for each $n\in\N$. Moreover, if $b\in \Q$, then system \eqref{eq:sist-LV-2}
	has a rational first integral and if $b \not \in  \Q$, then system \eqref{eq:sist-LV-2}
	has not a rational first integral.
\end{theorem}

{\bf Acknowledgements.}
The first author would like to thank  Universidad del B\'\i o-B\'\i o for the support provided during the master's studies at that institution, a period during which the results presented here were obtained.


\end{document}